\documentclass{ifacconf}

\usepackage{graphicx}      
\usepackage{natbib}        
\usepackage{amsmath}
\usepackage{amssymb}
\usepackage{xcolor}
\newcommand{\dz}{\operatorname{dz}}               
\newcommand{\sat}{\operatorname{sat}}             
\newcommand{\He}{\operatorname{He}}               
\usepackage{subfigure}
\usepackage{enumitem}

\begin{document}
	\begin{frontmatter}
		
		\title{Data-driven control of input saturated systems: a LMI-based approach} 
		
		\thanks[footnoteinfo]{This paper is partially supported by FAIR (Future Artificial Intelligence Research) project, funded by the NextGenerationEU program within the PNRR-PE-AI scheme (M4C2, Investment 1.3, Line on Artificial Intelligence), by the Italian Ministry of Enterprises and Made in Italy in the framework of the project 4DDS (4D Drone Swarms) under grant no. F/310097/01-04/X56 and by the PRIN PNRR project  P2022NB77E “A data-driven cooperative framework for the management of distributed energy and water resources” (CUP: D53D23016100001), funded by the NextGenerationEU program.}
		
		\author[First]{F. Porcari} 
		\author[Second]{V. Breschi} 
		\author[Third]{L. Zaccarian}
		\author[First]{S. Formentin}
		
		\address[First]{Dipartimento di Elettronica, Informazione e Bioingegneria, Politecnico di Milano, Milano, Italy (e-mail: federico.porcari@polimi.it)}
		\address[Second]{Department of Electrical Engineering, Eindhoven University of Technology, Eindhoven, The Netherlands}
		\address[Third]{LAAS-CNRS, University of Toulouse, Toulouse, France, and Dipartimento di Ingegneria Industriale, Università di Trento, Trento, Italy}
		
		\begin{abstract}                
This paper addresses three complex control challenges related to input-saturated systems from a data-driven perspective. Unlike the traditional two-stage process involving system identification and model-based control, the proposed approach eliminates the need for an explicit model description. The method combines data-based closed-loop representations, Lyapunov theory, instrumental variables, and a generalized sector condition to formulate data-driven linear matrix inequalities (LMIs). These LMIs are applied to maximize the origin's basin of attraction, minimize the closed-loop reachable set with bounded disturbances, and introduce a new data-driven $\ell_2$-gain minimization problem. Demonstrations on benchmark examples highlight the advantages and limitations of the proposed approach compared to an explicit identification of the system, emphasizing notable benefits in handling nonlinear dynamics.
		\end{abstract}
		
		
	\end{frontmatter}
	
	\section{Introduction}
The increasing availability of large datasets and the complexity of controlled systems have spurred interest in learning from data. This study explores the application of the data-driven paradigm to address input-saturated systems, a domain that has received limited attention despite the success of this paradigm in various control problems (see, e.g., \cite{Dor:23}). Existing model-reference, data-driven anti-windup approaches (\cite{Bre:20b}, \cite{Bre:20}) effectively handle input saturation but lack closed-loop stability guarantees. Recent contributions (\cite{Bre:23,Seu:23a,Seu:23b}) have attempted to translate established model-based design strategies (\cite{Tar:11}) into the data-driven realm for input-saturated systems, focusing on stabilizing state-feedback controllers.\\
In this study, a novel data-driven perspective is presented for three 
control algorithms tailored for input-saturated systems. The approach integrates a data-based closed-loop representation (\cite{Per:20}), Lyapunov theory, instrumental variables for noise handling (\cite{Bre:23}), and the generalized sector condition proposed in \cite{Sil:05}. Data-driven linear matrix inequalities (LMIs) are formulated, enabling the maximization of an estimate of the origin's basin of attraction, the minimization of the closed-loop reachable set with an energy-bounded input, and an alternative data-driven formulation for the nonlinear $\ell_2$-gain minimization problem 
already considered in \cite{Seu:23b}. The effectiveness of the strategies is demonstrated on an open-loop unstable benchmark (\cite{Bre:23}), comparing the closed-loop performance with an oracle controller designed using the true system dynamics. Emphasizing the maximization of the estimate of the basin of attraction, the study also discusses the performance of the proposed data-driven solution compared to the identification of the system dynamics on the same benchmark and on an attitude control problem for a quadcopter.\\
	The paper is organized as follows. The setup and the goal are introduced in Section~\ref{sec:problem_statement}, while the data-driven closed-loop representation is given in Section~\ref{sec:preliminaries}. The LMI-based data-driven strategies are then described in Section~\ref{sec:control_design}, while their effectiveness is assessed in Section~\ref{sec:numerical}. The comparison with model identification, both with and without modelling errors, is shown in Section~\ref{sec:result_comparison}, followed by some concluding remarks. \vspace{-.5cm}
	\paragraph*{Notation.} Let $\mathbb{R}$, $\mathbb{R}^{n}$ and $\mathbb{R}^{n \times m}$ be the set of real numbers, column vectors of length $n$ and $n \times m$ dimensional matrices, respectively. Given a vector $u \in \mathbb{R}^{n_{u}}$, we denote its $j$-th component as $u_{j}$, for $j=1,\ldots,n_{u}$, while we define the associated decentralized saturation function $\sat:\mathbb{R}^{n_u}\to\mathbb{R}^{n_u}$ as a function having components
			\begin{equation}\label{eq:saturation}
			\mathrm{sat}_j(u_{j})=\max\{\underline{u}_{j},\min\{\overline{u}_{j},u_{j}\}\},~~j=1,\ldots,n_{u},
		\end{equation}
	where $\overline{u}$ and $\underline{u}$ represent the upper and lower saturation bounds on $u$. Given $u \in \mathbb{R}^{n_{u}}$, we further define the dead-zone nonlinearity as:
	\begin{equation}\label{eq:dead_zone}
		\mathrm{dz}(u)=u-\mathrm{sat}(u).
	\end{equation}
	Given a positive definite matrix $Q \in \mathbb{R}^{n_x \times n_x}$ and a scalar $s > 0$, the ellipsoidal subset $\mathcal{E}(Q,s)$ is defined as 
	\begin{equation}\label{eq:ellipsoid}
		\mathcal{E}(Q,s) = \left\{ x:x^\top Q^{-1} x \leq s^2 \right\}.
	\end{equation} 
	The shorthand notation $x^+ = Ax$ is used instead of $x(k+1) = A x(k)$. 
	For a full-rank matrix $Z\in \mathbb{R}^{n \times m}$, $Z^\dagger$ denotes its right pseudo-inverse while $Z_{i}$ indicates its $i$-th row, with $i=1,\ldots,n$. For a square matrix $A$, $c(A)$ denotes its condition number and  $\He(A)$ is twice the symmetric part of $A$, i.e., $\He(A) = A + A^\top$. For a discrete-time signal $\xi(k) \in \mathbb{R}^{n_{\xi}}$ and for any $k_0, k_1, L \in \mathbb{N}$ such that $0 \leq k_0 < k_1$ and $L \leq k_1 - k_0 + 1$, the Hankel matrix $\Xi_{k_0, L, k_1} \in \mathbb{R}^{n_{\xi}L\times k_{1}-k_{0}-L+1}$ is defined as 
	\begin{equation}
		\Xi_{k_0, L, k_1} \!=\!\! \begin{bmatrix}
			\xi(k_0) & \xi(k_0\!+\!1) & \ldots & \xi(k_1\!-\!L\!+\!1) \\
			\xi(k_0\!+\!1) & \xi(k_0\!+\!2) & \ldots & \xi(k_1\!-\!L\!+\!2) \\
			\vdots & \vdots & \ddots & \vdots \\
			\xi(k_0\!+\!L\!-\!1) & \xi(k_0\!+\!L) & \ldots & \xi(k_1)
		\end{bmatrix}\!,
	\end{equation}
	while $\Xi_{k_0,k_1} = \Xi_{k_0,1,k_1}$ is the single-row Hankel matrix.
	
	\section{Setting and control objectives}\label{sec:problem_statement}
	Consider a discrete-time, linear time invariant system characterized by the state and output equations
	\begin{subequations}
		\begin{equation}\label{eq:plant_equations}
        \begin{aligned}
			&x^{+}=Ax+B\mathrm{sat}(u)+w = Ax + Bv + w,\\
			&y=x+e,
        \end{aligned}
		\end{equation}
	where $x \in \mathbb{R}^{n_x}$ is its state, which is assumed to be \emph{fully measurable} but corrupted by a zero-mean measurement white noise $e \in \mathbb{R}^{n_x}$ with  positive definite covariance $\Sigma_{e} \in \mathbb{R}^{n_x \times n_x}$, and $w \in \mathbb{R}^{n_x}$ is an uncontrollable but known disturbance 
	with \emph{bounded energy}, i.e.,
		\begin{equation}\label{eq:ell2_bound}
		\|w\|_{2} \leq s,~~s\geq 0.
	\end{equation}  
	Let $u \in \mathbb{R}^{n_u}$ be the system's input and $v=\mathrm{sat}(u)$ be defined as in \eqref{eq:saturation}, with the saturation here assumed to be symmetric with respect to the origin, namely $\overline{u}_{j}=\underline{u}_{j}$ for all $j \in \{1,\ldots,n_{u}\}$, and the bound $\overline{u}$ being \emph{known}.\\
	\end{subequations}
	Let us suppose that the system's performance is encoded into the 
	signal $z \in \mathbb{R}^{n_{z}}$ defined as:
	\begin{equation}\label{eq:z}
		z= Cx+D_{w}w+D_{u}v = Cx+D_{w}w+D_{u}\mathrm{sat}(u).
	\end{equation}
	Our goal is to design a static state-feedback controller
	\begin{equation}\label{eq:feedback}
		u=Kx,~~K \in \mathbb{R}^{n_{u}\times n_{x}},
	\end{equation}
	such that the origin of the resulting closed-loop system is \emph{asymptotically stable}, while one of the following is satisfied:
	\begin{itemize}[leftmargin=.45in]
		\item[(BoA)] for $w=0$, an ellipsoidal estimate $\mathcal{E}(Q,1)$ (see \eqref{eq:ellipsoid}) of the origin's basin of attraction is maximized, and the closed-loop response satisfies
		\begin{equation}\label{eq:convergence}
		\lvert x(t) \rvert < \eta^{t} \sqrt{c(Q)} \, \lvert x(0) \rvert ,~~\eta \in (0,1], \forall t \geq 0;
		\end{equation}
		\item[(RS)] an ellipsoidal estimate $\mathcal{E}(Q,s)$ of the closed-loop reachable set $\mathcal{S}(s) \subset \mathbb{R}^{n_x}$ defined as
		\begin{equation}
			\mathcal{S}(s) \!=\! \{x\!:\! x(0)\!=\!0 \!\Rightarrow\! x(t) \!\in\! \mathcal{S}, \forall t \!\geq\! 0, \forall w\!:\! \|w\|_{2}\!\leq\! s\},
		\end{equation}
		 is minimized for a given bound $s$ on $\|w\|_2$;
		\item[($\ell_2$)] an estimate $\gamma(s)$ of the closed-loop $\ell_{2}$-gain from $w$ to $z$, defined for $x(0) = 0$ as
        \begin{equation}
            \sup_{\substack{w \neq 0 \\ \lVert w \rVert_2 \leq s}} \frac{\lVert z \rVert_2}{\lVert w \rVert_2},
        \end{equation}
        is minimized for a given bound $s$ on $\|w\|_2$.
	\end{itemize}
	
	While $C \in \mathbb{R}^{n_z \times n_{x}}$, $D_{u} \in \mathbb{R}^{n_z \times n_{u}}$ and $D_{w} \in \mathbb{R}^{n_z \times n_{x}}$ are fixed based on the control objective, assume that a model for the controlled system is \emph{not known}. Nonetheless, assume we have access to a dataset
	\begin{equation}\label{eq:dataset}
		\mathcal{D}_{T}=
		\{w^{d}(k),v^{d}(k),y^{d}(k)\}_{k=0}^{T},
	\end{equation}
	of length
	\begin{equation}\label{eq:length}
		T \geq (n_{u}+1)n_{x}+n_{u},
	\end{equation}
	obtained by feeding the plant with a saturated input designed to be \emph{persistently exciting} of order $n_{x}+1$, where persistence of excitation is defined as follows.
	\begin{defn}[Persistence of excitation]
		A discrete time signal $v(k) \in \mathbb{R}^{n_u}$, $k=0,\ldots,T$ is persistently exciting of order $L$ if the Hankel matrix $V_{0,L,T-1}$ is full row rank. 
	\end{defn}
	This implies that the input $u^d$ is designed to be persistently exciting and that the rank property is not lost in $v^d$ due to the saturation effect. In this setting, our goal translates into designing from data a stabilizing state feedback gain $K$ as in \eqref{eq:feedback}, satisfying one of the aforementioned design conditions: an endeavor that we aim here to undertake via a \textquotedblleft hybrid\textquotedblright \ approach stemming from the one presented in \cite{Bre:23}.

	\section{Closed-loop data-driven description}\label{sec:preliminaries}
	Following the footsteps and tricks already adopted in \cite{Per:20,Bre:23} we retrieve a data-driven description of the system. To this end, let us assume that we have collected a second set of outputs $\{\tilde{y}^{d}(k)\}_{k=0}^{T}$ by feeding the plant with the same input sequence used to gather $\{y^{d}(k)\}_{k=0}^{T}$. For both collected datasets, the persistence of excitation condition on $v^{d}$ and the imposed lower-bound on $T$ (see \eqref{eq:length}) imply that the following holds~(see \cite{Wil:05}):  
	\begin{equation} \label{eq:rank_condition} 
		\textrm{rank} \left( \begin{bmatrix}
			V_{0,T-1}^{d} \\ Y_{0,T-1}^{d}
		\end{bmatrix} \right) = 
        \textrm{rank} \left( \begin{bmatrix}
			V_{0,T-1}^{d} \\ \tilde{Y}_{0,T-1}^{d}
		\end{bmatrix} \right) = n_u + n_x.
	\end{equation}    
	Moreover, the availability of the second dataset allows us to construct an \emph{instrumental} Hankel matrix $Z_{0,T-1}^{d}$ as
	\begin{equation}\label{eq:instrumental_variable}
		Z^{d}_{0,T-1} = \begin{bmatrix}
			V_{0,T-1}^{d} \\ \tilde{Y}_{0,T-1}^{d}
		\end{bmatrix},
	\end{equation}
    which satisfies the rank condition \eqref{eq:rank_condition}.
    By leveraging the non-correlation property of white noises, i.e. 
    \begin{equation}\label{eq:correlation}
        \mathbb{E}\left\{ e(t) \tilde{e}(\tau) \right\} = 0, \qquad \forall t, \tau,
    \end{equation} 
    where $e$ and $\tilde{e}$ are, respectively, the noises acting on $y^d$ and on $\tilde{y}^d$, the matrix \eqref{eq:instrumental_variable} can be used to extract an (asymptotically unbiased) estimate of the open-loop system's matrices in \eqref{eq:plant_equations} via the instrumental variable least-squares problem (\cite{Sod:02})
	\begin{align}\label{eq:ol_est}
		\nonumber \begin{bmatrix}
			\hat{B} & \hat{A}
		\end{bmatrix}&=\underset{\begin{bmatrix} B & A\end{bmatrix}}{\arg\!\min}~\bigg\|\left(\bar{Y}_{1,T}^{d}-\begin{bmatrix} B & A\end{bmatrix}\begin{bmatrix}
			V_{0,T-1}^{d}\\
			Y_{0,T-1}^{d}
		\end{bmatrix}\right)(Z^{d}_{0,T-1})^{\top}\bigg\|_{2}^{2}\\
		& =\bar{Y}_{1,T}^{d}(Z^{d}_{0,T-1})^{\top}\left(\begin{bmatrix}
		V_{0,T-1}^{d}\\
		Y_{0,T-1}^{d}
	\end{bmatrix}(Z^{d}_{0,T-1})^{\!\top}\right)^{\!-1}\!\!,
	\end{align}
	where $\bar{Y}_{1,T}^{d}=Y_{1,T}^{d}-W_{0,T-1}^{d}$ is known because $w$ in \eqref{eq:plant_equations} is known and the inverse of the data matrices is well-defined thanks to the full rank condition in \eqref{eq:rank_condition}. In turn, this leads to the data-based representation
	\begin{equation}\label{eq:DD_ol}
		x^{+}=\bar{Y}_{1,T}^{d}(Z^{d}_{0,T-1})^{\top}\!\left(\begin{bmatrix}
			V_{0,T-1}^{d}\\
			Y_{0,T-1}^{d}
		\end{bmatrix}(Z^{d}_{0,T-1})^{\!\top}\!\!\right)^{\!\!-1}\begin{bmatrix}
		\mathrm{sat}(u)\\
		x
	\end{bmatrix}\!+w,
	\end{equation}
	as formalized in the following lemma.
	\begin{lem}[Open-loop description]\label{lem:open_loop_description}
			Given $\mathcal{D}_{T}$ in \eqref{eq:dataset} and a fresh set of outputs $\{\tilde{y}^{d}(k)\}_{k=0}^{T}$ collected by exploiting the same input used to construct $\mathcal{D}_{T}$, the dynamics of the system with saturated inputs \eqref{eq:plant_equations} can be expressed as in \eqref{eq:DD_ol} with an error that vanishes as the length of the dataset increases, namely
			\begin{equation}\label{eq:unbiased}
				\lim_{T \rightarrow \infty} \begin{bmatrix}
					\hat{B} & \hat{A}
				\end{bmatrix} \overset{\mathrm{a.s.}}{\longrightarrow} \begin{bmatrix}
				B & A
			\end{bmatrix}.
			\end{equation}
	\end{lem}
	\begin{pf}
		To show that \eqref{eq:unbiased} holds, we can follow the same steps carried out in \citep[Proof of Lemma 4]{Bre:23}. Hence, let us replace $\bar{Y}_{1,T}^{d}$ with
		\begin{equation*}
			\bar{Y}_{1,T}^{d}=\begin{bmatrix}
				B & A
			\end{bmatrix}\begin{bmatrix}
			V_{0,T-1}^{d}\\
			Y_{0,T-1}^{d}
		\end{bmatrix}+\underbrace{E_{1,T}-AE_{0,T-1}}_{\Delta E_{0,T-1}},
		\end{equation*}
		where $E_{1,T}$ and $E_{0,T-1}$ are the Hankel matrices associated with noise $e$. Accordingly, \eqref{eq:ol_est} can be equivalently rewritten as
		\begin{equation*}
			\begin{bmatrix}
				\hat{B} & \hat{A}
			\end{bmatrix}\!=\!\begin{bmatrix}
			B & A
		\end{bmatrix}\!+\!\Delta E_{0,T\!-1}(Z^{d}_{0,T\!-1})^{\!\top}\!\left(\begin{bmatrix}
		V_{0,T-1}^{d}\\
		Y_{0,T-1}^{d}
	\end{bmatrix}\!(Z^{d}_{0,T\!-1})^{\!\top}\!\right)^{\!\!-1}\!\!\!\!.
		\end{equation*}
	In turn, considering $\tilde{Y}_{0,T-1}^d=X_{0,T-1}+\tilde{E}_{0,T-1}$ and that as $T \to \infty$ the sample mean becomes the expected value, due to the fact that $X_{0,T-1}$ and $V_{0,T-1}^d$ do not depend on the measurement noise and, thus, the expected value of their product with $\Delta E_{0,T-1}$ is zero, from \eqref{eq:correlation} it follows that $\Delta E_{0,T-1}(Z^{d}_{0,T-1})^{\top}$ in the second term of the former equality satisfies
	\begin{equation*}
		\frac{\Delta E_{0,T\!-1}(Z^{d}_{0,T\!-1})^\top}{T} \underset{T \rightarrow \infty}{\longrightarrow} \mathbb{E}\left\{\Delta E_{0,T\!-1}\begin{bmatrix}
			0 \\ \tilde{E}_{0,T-1}
		\end{bmatrix}^{\top}\right\}=0,
	\end{equation*}
	thus concluding the proof. \hfill $\blacksquare$
	\end{pf}
	
	Based on Lemma \ref{lem:open_loop_description}, the closed-loop dynamics \eqref{eq:plant_equations} and \eqref{eq:feedback} can be written as follows:
	\begin{subequations}\label{eq:data_based_close_loop}
	\begin{equation}
		x^+ = A_{\mathrm{cl}}^{d}(G) x - B_{\mathrm{cl}}^{d} \dz(u) + w,
	\end{equation}
	where 
	\begin{align}
			& A_{\mathrm{cl}}^{d}(G) = \bar{Y}_{1,T}^{d} (Z^{d}_{0,T-1})^\top G, \\
			& B_{\mathrm{cl}}^{d} \!=\! \bar{Y}_{1,T}^{d} (Z^{d}_{0,T-1})^{\!\top}\! \left( \left[ \begin{array}{c}
				V_{0,T-1}^{d} \\ Y_{0,T-1}^{d}
			\end{array} \right] (Z^{d}_{0,T-1})^{\!\top} \right)^{\!\!-1} \begin{bmatrix}
				I \\ 0
			\end{bmatrix},
	\end{align}
	with $G \in \mathbb{R}^{(n_u+n_x) \times n_x}$ satisfying
	\begin{equation}\label{eq:consistency_origin}
		\begin{bmatrix}
			K \\ I
		\end{bmatrix} = \begin{bmatrix}
			V_{0,T-1}^{d} \\ Y_{0,T-1}^{d}
		\end{bmatrix} (Z^{d}_{0,T-1})^\top G.
	\end{equation}
	\end{subequations}
	This result is formally stated in the following lemma.
		\begin{lem}[Closed-loop description]\label{lemma2}
		Given $\mathcal{D}_{T}$ in \eqref{eq:dataset} and a set of outputs $\{\tilde{y}^{d}(k)\}_{k=0}^{T}$ collected by exploiting the same input used to construct $\mathcal{D}_T$, \eqref{eq:data_based_close_loop} asymptotically (i.e., as $T \rightarrow \infty$) coincides with the closed-loop interconnection \eqref{eq:plant_equations}, \eqref{eq:feedback}. 
	\end{lem}
	\begin{pf}
		The proof easily follows similar steps to those in \cite[Theorem 3]{Bre:23} and is omitted. It hinges upon the fact that \eqref{eq:ol_est} and  \eqref{eq:consistency_origin} provide 
        \begin{align*}
           \begin{bmatrix}
			\hat{B} & \hat{A}
		\end{bmatrix} \begin{bmatrix}
			K \\ I
    		\end{bmatrix} = & \; \bar{Y}_{1,T}^{d}(Z^{d}_{0,T-1})^{\!\top}\left(\begin{bmatrix}
    		V_{0,T-1}^{d}\\
    		Y_{0,T-1}^{d}
	    \end{bmatrix}(Z^{d}_{0,T-1})^{\!\top}\right)^{\!-1}\!\! \\
        & \times \begin{bmatrix}
			V_{0,T-1}^{d} \\ Y_{0,T-1}^{d}
		\end{bmatrix} (Z^{d}_{0,T-1})^\top G.
        \end{align*}\hfill $\blacksquare$ 
	\end{pf}
	Note that $G$ in \eqref{eq:data_based_close_loop} plays the role of a tuning parameter that uniquely defines $K$ through \eqref{eq:consistency_origin}. While $A_{\mathrm{cl}}^{d}$ in \eqref{eq:data_based_close_loop} is fully characterized from data, we still need an estimate of $B$, i.e., $B_{cl}^{d}$, to describe the effect of the input saturation.
	
	\begin{rem}
		The products 
		\begin{equation*}
			\bar{Y}_{1,T}^{d} \left(Z^{d}_{0,T-1} \right)^\top,~~~~~V_{0,T-1}^{d} \left( Z^{d}_{0,T-1} \right)^\top,
		\end{equation*}
		may be large in magnitude when the dataset is large. To avoid numerical problems, it is thus advisable to use the normalized instrumental variable
		\begin{equation}
			\bar{Z}_{0,T-1}^{d} = \frac{1}{T} \begin{bmatrix}
				V_{0,T-1}^{d} \\ \tilde{Y}_{0,T-1}^{d}
			\end{bmatrix}.
		\end{equation}
		As it can be easily proven by replacing $Z^{d}_{0,T-1}$ with $\bar{Z}_{0,T-1}$ in \eqref{eq:ol_est}, using the normalized instrument still leads to a closed-loop representation asymptotically equivalent to the one in \eqref{eq:data_based_close_loop}.
	\end{rem}
	
	\section{Data-driven design strategies}\label{sec:control_design}
    \subsection{Regional quadratic certificates}
	By relying on the data-based closed-loop representation in \eqref{eq:data_based_close_loop} we can now solve the design problems listed in Section~\ref{sec:problem_statement}. To this end, let us introduce the candidate Lyapunov function
	\begin{equation}\label{eq:lypaunov}
		V(x)=x^{\top}Px,
	\end{equation}
	with $P \in \mathbb{R}^{n_x \times n_x}$ being symmetric and positive definite. Then, let us recall that global asymptotic stability in the presence of saturation can only be achieved with non-exponentially unstable plants (\cite{Las:93}). This condition is rather limiting since the data-driven closed-loop description in \eqref{eq:data_based_close_loop} can also be retrieved for exponentially unstable plants, by carrying out suitably designed closed-loop experiments\footnote{This would require a pre-designed stabilizing controller and it would make noise handling more challenging. Therefore, the treatment of this case is left for future works.}. In designing the feedback gain in \eqref{eq:feedback} we thus look for regional stability, by considering the \emph{generalized sector condition} (\cite{Sil:05}) stating that for any $u \in \mathbb{R}^{n_u}$, $H \in \mathbb{R}^{n_u \times n_x}$ and diagonal, positive definite $W \in \mathbb{R}^{n_u \times n_u}$, the following holds 
	\begin{equation}\label{eq:generalized_sector}
		\dz(Hx)=0 \Rightarrow \dz(u)^\top W (u - \dz(u) + Hx) \geq 0.
	\end{equation}
	The generalized sector condition \eqref{eq:generalized_sector} can be exploited in an ellipsoidal estimate $\mathcal{E}(Q,s)$, as long as $\mathcal{E}(Q,s)$ is contained in the subset of $\mathbb{R}^{n_x}$ where $\dz(Hx)=0$. This set inclusion can be enforced by imposing 
	\begin{equation}\label{eq:ellipsoid_in_deadzone}
		x^\top \frac{H_j^\top H_j}{\overline{u}_j^2} x < \frac{1}{s^2}x^\top Px,~~\forall j \in \{1,\ldots,n_u\}, 
	\end{equation}
    which implies $x^\top P x \leq s^2 \Rightarrow \lvert H_j x \rvert^2 \leq \overline{u}_j^2$ for each $j$.
    By a Schur complement and a congruence transformation, condition \eqref{eq:ellipsoid_in_deadzone} can be rewritten as the following set of LMIs in the features of the dead-zone, the saturation and the candidate Lyapunov function \eqref{eq:lypaunov}:
	\begin{equation}
		\label{eq:saturation_LMI}
		\left[ \begin{array}{cc}
			Q & \;\;N_j^\top \\[2mm]
			N_j & \;\;\overline{u}_j^2 / s^2
		\end{array} \right] \succ 0, \qquad \forall \, j \in \left\{ 1, \ldots, n_u \right\},
	\end{equation}
	which is linear in the transformed decision variables $N = HQ \in \mathbb{R}^{n_u \times n_x}$ and $Q=P^{-1}$. 
	\subsection{BoA estimate with certified convergence rate}
	Let us initially assume that no external signal affects the system except for the controlled input $u$ (i.e., $w=0$). To find an ellipsoidal estimate of the origin's basin of attraction while guaranteeing a desired local exponential convergence rate for the closed-loop solutions (see \eqref{eq:convergence}) one has to concurrently solve \eqref{eq:saturation_LMI} for $s=1$ (in fact any value of $s$ should work up to a rescaling of $Q$) jointly with
	\begin{equation}\label{eq:DoA_LMI}
		 \He \begin{bmatrix}
			-\tfrac{\eta}{2} Q & 0 & 0\\
			V_{0,T-1}^{d} \left( Z^{d}_{0,T-1} \right)^\top F + N & - M & 0 \\
			\bar{Y}_{1,T}^{d} \left( Z^{d}_{0,T-1} \right)^\top F & -B_{\mathrm{cl}}^{d} M & -\tfrac{\eta}{2} Q
		\end{bmatrix} \prec 0,
	\end{equation}  
	where $M=W^{-1} \in \mathbb{R}^{n_u \times n_u}$ is any diagonal and positive definite matrix issued from \eqref{eq:generalized_sector}, and $F \in \mathbb{R}^{(n_u+n_x)\times n_x}$ satisfies the consistency condition
	\begin{equation}\label{eq:consistency_equation}
		Y_{0,T-1}^{d} \left(Z^{d}_{0,T-1} \right)^\top F = Q,
	\end{equation}
	according to the second block component of \eqref{eq:consistency_origin}. This leads in the following design result.
	\begin{thm}\label{thm1}
		Given \eqref{eq:plant_equations}, a dataset $\mathcal{D}_{T}$ satisfying \eqref{eq:length} and \eqref{eq:rank_condition}, and the instrument $Z^{d}_{0,T-1}$ in \eqref{eq:instrumental_variable}, if there exist matrices $Q=Q^{\top} \in \mathbb{R}^{n_x \times n_x}$, $F \in \mathbb{R}^{(n_{u}+n_{x})\times n_{x}}$, $N \in \mathbb{R}^{n_{u}\times n_{x}}$ and $M \in \mathbb{R}^{n_u\times n_u}$ diagonal satisfying \eqref{eq:saturation_LMI} for $s=1$, \eqref{eq:DoA_LMI} and \eqref{eq:consistency_equation}, then the feedback gain
		\begin{equation}\label{eq:data_driven_K}
			K = V_{0,T-1}^{d} \left( Z^{d}_{0,T-1} \right)^{\top} F Q^{-1},
		\end{equation} 
		 asymptotically (i.e., for $T \rightarrow \infty$) guarantees the following: $(i)$ closed-loop exponential stability of the origin with rate $\eta$, namely solutions satisfy \eqref{eq:convergence} in the basin of attraction of the origin; $(ii)$ $\mathcal{E}(Q,1)$ is contained in the origin's basin of attraction.
	\end{thm}
	\begin{pf}
	Consider the Lyapunov candidate function in \eqref{eq:lypaunov}, with $P=P^{\top}=Q^{-1}$ being positive definite thanks to \eqref{eq:DoA_LMI} and Sylvester's criterion. We study the exponential stability of the origin by evaluating the sign of
	\begin{equation}\label{eq:delta_Veta}
		\Delta_{\eta} V(x)=(x^{+})^{\top}Px^{+}-\eta^{2}x^{\top}Px \leq 0.
	\end{equation}
    The proof then unfolds as the one of Theorem~1 in \cite{Bre:23}. Specifically, starting from $Q$, $F$, $M$ and $N$ satisfying \eqref{eq:saturation_LMI}-\eqref{eq:DoA_LMI} as assumed in the statement of the theorem, select $G=FQ^{-1}$ so that, from \eqref{eq:consistency_equation} and \eqref{eq:data_driven_K}, we obtain that \eqref{eq:consistency_origin} holds for this specific selection of $G$ and $K$. Then, according to Lemma~\ref{lemma2}, \eqref{eq:delta_Veta} is asymptotically equivalent to
    \begin{equation} \label{eq:expressionVdeltaeta}
	\begin{aligned}
		\eta^{-1} \Delta_{\eta} V(x)= & \\ 
            \xi^{\!\top}\left(\begin{bmatrix}
			(A_{\mathrm{cl}}^{d}(G))^{\!\top}\\
			-(B_{\mathrm{cl}}^{d})^{\!\top}
		\end{bmatrix}\! \right.&\left.\!\frac{P}{\eta}\!\begin{bmatrix}
			A_{\mathrm{cl}}^{d}(G) &
			-B_{\mathrm{cl}}^{d}
	\end{bmatrix}\!-\!\begin{bmatrix}
		\eta P & 0\\
		0 & 0
	\end{bmatrix} \right)\xi \leq 0,
	\end{aligned}
    \end{equation}
	where $\xi\!=\!\begin{bmatrix}
		x\\
		\dz(u)
	\end{bmatrix}$. In turn, for any $x$ satisfying $\dz(Hx)\!=\!0$, \eqref{eq:expressionVdeltaeta} can be upper-bounded by using the sector condition \eqref{eq:generalized_sector} as follows:
		\begin{equation*}
			\eta^{-1} \Delta_{\eta} V(x) \leq \xi^{\top}\Xi_{\eta}\xi,
		\end{equation*} 
	where
	\begin{equation}\label{eq:Xi_eta}
		\Xi_{\eta}=\mathrm{He}\begin{bmatrix}
			-\frac{\eta}{2}P & 0 & 0\\
			WK+WH & -W & 0\\
				\bar{Y}_{1,T}^{d} \left( Z^{d}_{0,T-1} \right)^\top G \;\; & -B_{\mathrm{cl}}^{d} \; & -\frac{\eta}{2}P^{-1}
		\end{bmatrix}.
	\end{equation}
	Selecting the variables appearing in \eqref{eq:Xi_eta} as
	\begin{align*}
		& P=Q^{-1}>0,~~G=FQ^{-1}, \\
		&W=M^{-1}>0 \mbox{ diagonal},~~H=NQ^{-1},
	\end{align*}
	using Schur's complements and pre- and post- multiplying the resulting matrix by $\mathrm{diag}\{P^{-1},W^{-1},I\}$, we obtain from \eqref{eq:DoA_LMI} that $\Xi_{\eta}<0$, thus showing that
	\begin{equation*}
		\dz(Hz)=0 \Rightarrow \eta^{-1} \Delta_{\eta} V(x) \leq \lambda_{max}(\Xi_{\eta})|x|^{2}<0,
	\end{equation*} 
	with $\lambda_{max}(\Xi_{\eta})<0$ denoting the largest eigenvalue of (the negative definite matrix) $\Xi_{\eta}$. Performing the same manipulations carried out in \cite[Proof of Theorem 1]{Bre:23}, it finally follows that $x \in \mathcal{E}(Q,1)\setminus \{0\} \Rightarrow \Delta_{\eta} V(x)<0$. This further implies that
	\begin{equation}\label{eq:time_iter}
		V(x^{+})\leq \eta^{2}V(x) \Rightarrow V(x(t))\leq \eta^{2t}V(x(0)),
	\end{equation}
	where the inequality on the right-hand side of the implication results from the propagation over time of the condition on the left-hand side. Based on the chosen Lyapunov candidate it further holds that
	\begin{equation*}
		\lambda_{\mathrm{min}}(P)|x|^{2}\leq V(x)\leq \lambda_{\mathrm{max}}(P)|x|^{2},
	\end{equation*}
	where $\lambda_{\mathrm{min}}(P),\lambda_{\mathrm{max}}(P)>0$ are the minimum and maximum eigenvalues of $P$, respectively. Combining the right-hand side of \eqref{eq:time_iter} with the previous bounds, we obtain
	\begin{equation*}
		\lambda_{\mathrm{min}}(P)|x(t)|^{2}\leq \eta^{2t} \lambda_{\mathrm{max}}(P)|x(0)|^{2},
	\end{equation*}
	entailing that the following also holds:
	\begin{align*}
			|x(t)|^{2}\leq \eta^{2t} \frac{\lambda_{\mathrm{max}}(P)}{\lambda_{\mathrm{min}}(P)}|x(0)|^{2}&=\eta^{2t}c(P)|x(0)|^{2}\\
			&=\eta^{2t}c(Q)|x(0)|^{2},
	\end{align*}
	where $c(P)$ is the condition number of $P$ and $c(P)=c(Q)$ because $P=Q^{-1}$. The closed-loop convergence rate in \eqref{eq:convergence} is then straightforwardly proven by taking the square root on both sides of the previous inequality, thus concluding the proof. \hfill $\blacksquare$
	\end{pf}
	As a direct translation from the model-based context (see, e.g., {\cite{Tar:11}}, Chapter~2), the first of our design problems can thus be tackled in a data-based fashion by solving the following Semidefinite Program (SDP):
 	\begin{equation}\label{eq:DoA_optimization}
 	\begin{aligned}
 		& \underset{\alpha, Q, U, F, N}{\mathrm{minimize}}  ~~ -\alpha \\
 		&\quad \quad~ \mathrm{s.t.}  \quad \alpha I \preceq Q, \\
 		& ~ \qquad \qquad \eqref{eq:saturation_LMI}, \eqref{eq:DoA_LMI},\eqref{eq:consistency_equation}.
 	\end{aligned}
	 \end{equation}
	Note that none of the variables to be optimized has a size that depends on $T$, thanks to the use of the instrument in \eqref{eq:instrumental_variable}. Therefore, $\mathcal{D}_{T}$ can in principle be arbitrarily long, which is desirable to attain the equivalence between the model-based and data-driven system's description.  
	
	\subsection{Minimized reachable set from bounded $w$}
	Let us now study the effect of the exogenous signal $w$ in \eqref{eq:plant_equations}. In this setting, the closed-loop stability of the origin and the fact that $\mathcal{E}(Q,s)$ is an outer approximation of the reachable set are asymptotically guaranteed by imposing:
		\begin{equation}\label{eq:reachable_set_LMI}
		\He\! \begin{bmatrix}
			-\tfrac{1}{2} Q \!&\! 0 & 0 & 0 \\
			V_{0,T-1}^{d}\! \left(Z^{d}_{0,T-1} \right)^{\!\top}\! \!F \!+\! N \!&\! - M & 0 & 0 \\
			0 \!&\! 0 & -\tfrac{1}{2} I & 0 \\
			\bar{Y}_{1,T}^{d} \!\left(Z^{d}_{0,T-1} \right)^{\!\top}\! \!F \!&\! -B_{\mathrm{cl}}^{d} M & I & -\tfrac{1}{2} Q
		\end{bmatrix} \!\!\prec\! 0,
	\end{equation}
	along with \eqref{eq:consistency_equation}, leading to the following design result.
	\begin{thm}\label{thm2}
		Let $w$ in \eqref{eq:plant_equations} verify \eqref{eq:ell2_bound}. Given the dataset $\mathcal{D}_{T}$ satisfying \eqref{eq:length} and \eqref{eq:rank_condition}, if there exist $Q=Q^{\top}$, a diagonal matrix $M$, $F \in \mathbb{R}^{(n_{u}+n_{x})\times n_{x}}$ and $N \in \mathbb{R}^{n_{u}\times n_{x}}$ satisfying \eqref{eq:saturation_LMI},  \eqref{eq:consistency_equation} and \eqref{eq:reachable_set_LMI}, then the feedback gain $K$ in \eqref{eq:data_driven_K} asymptotically (i.e., for $T \rightarrow \infty$) guarantees $(i)$ closed-loop exponential stability of the origin, and that $(ii)$ $\mathcal{E}(Q,s)$ is an outer approximation of the closed-loop system's reachable set $\mathcal{S}$ from $x(0)=0$ and for any $w$ satisfying \eqref{eq:ell2_bound}.
	\end{thm}
	\begin{pf}
		Based on the element in position $(1,1)$ of \eqref{eq:reachable_set_LMI}, the Lyapunov candidate function in \eqref{eq:lypaunov}, with $P=Q^{-1}$ is positive definite. By leveraging the asymptotic equivalence between the \textquotedblleft true\textquotedblright \ closed-loop system and its data-based description in \eqref{eq:data_based_close_loop}, we can characterize $\Delta V(x,w) = V(x^+) - V(x)$ as follows:
		\begin{equation}\label{eq:equality_proof2}
			\Delta V(x,w)=\tilde{\xi}^{\top}\Gamma\tilde{\xi},
		\end{equation}
		where 
		\begin{equation*}
        \begin{aligned}
			\Gamma &= \begin{bmatrix}
				(A_{\mathrm{cl}}^{d}(G))^{\!\top}\\
				-(B_{\mathrm{cl}}^{d})^{\!\top}\\
				I
			\end{bmatrix} P \begin{bmatrix}
				A_{\mathrm{cl}}^{d}(G) &
				-B_{\mathrm{cl}}^{d} & I
			\end{bmatrix} - \begin{bmatrix}
				P & 0 & 0\\
				0 & 0 & 0\\
				0 & 0 & 0
			\end{bmatrix},
            \\\tilde{\xi} &= \begin{bmatrix}
				x\\
				\dz(u)\\
				w
		\end{bmatrix}.
        \end{aligned}
		\end{equation*}
		With $w=0$, since \eqref{eq:reachable_set_LMI} implies \eqref{eq:DoA_LMI}, then exponential stability follows from Theorem~\ref{thm1}, whose steps imply $\Delta V(x,0) \leq -(1-\eta^2)V(x)$. For $w\neq0$, the same steps as those of Theorem~\ref{thm1}, together with \eqref{eq:equality_proof2} and inequality \eqref{eq:reachable_set_LMI}, imply
        \begin{equation}\label{eq:delta_V}
            \Delta V(x,w) \leq -(1-\eta^2)V(x) + w^\top w \leq w^\top w.
        \end{equation}
        Summing \eqref{eq:delta_V} and using $x(0)=0 \Rightarrow V(x(0))=0$, we get 
        \begin{align*}
		      V(x(t)) & = V(x(t)) - V(x(0)) \\ & = 
            \sum_{\tau=0}^{t-1} V(x(\tau+1)) - V(x(\tau)) = \sum_{\tau=0}^{t-1} \Delta V(x(\tau)) \\ & \leq \sum_{\tau=0}^{t-1} \lvert w(\tau) \rvert^2 \leq s^2,
        \end{align*}
		implying that enforcing \eqref{eq:saturation_LMI} and \eqref{eq:reachable_set_LMI} asymptotically guarantees that $x(t) \in \mathcal{E}(Q,s)$ since $Q=P^{-1}$.
	\hfill $\blacksquare$
	\end{pf}
	Note that the expression of the feedback gain in this design result corresponds to the one in \eqref{eq:data_driven_K}, but its actual value is shaped by matrices that verify the new LMI in \eqref{eq:reachable_set_LMI}. In turn, this result allows us to formulate the following data-driven SDP for the reachable set minimization:
	\begin{equation}\label{eq:trace_min}
		\begin{aligned}
			& \underset{Q, U, F, N}{\mathrm{minimize}}  ~~ \textrm{trace}(Q) \\
			&\quad \quad~  \mathrm{s.t.}  \quad \eqref{eq:saturation_LMI}, \eqref{eq:consistency_equation}, \eqref{eq:reachable_set_LMI},
		\end{aligned}
	\end{equation}
    which is computationally tractable even for large datasets, thanks to the introduction of the instrument.
 
	\subsection{Minimum $\ell_2$-gain}
	By still assuming the bound in \eqref{eq:ell2_bound}, let us now consider the problem of minimizing the $\ell_{2}$-gain $\gamma(s)$ of the closed-loop system from the exogenous input $w$ to the performance variable $z$, namely
	\begin{equation}\label{eq:l2_gain_condition}
		x(0) = 0 \Rightarrow \| z \|_2 \leq \gamma(s) \| w \|_2 ,~ \quad \forall w: \lVert w \rVert_2 \leq s.
	\end{equation} 
 \begin{table*}[!h]
		{\normalsize
			\begin{equation} \label{eq:l2_gain_LMI}
				\He\begin{bmatrix}
					-\tfrac{1}{2} Q & 0 & 0 & 0 & 0\\
					V_{0,T-1}^{d} \left( Z^{d}_{0,T-1} \right)^{\!\top}\!\! F + N & - M & 0 & 0 & 0 \\
					0 & 0 & -\tfrac{1}{2} I & 0 & 0\\
					C Q + D_u V_{0,T-1}^{d} \left( Z^{d}_{0,T-1} \right)^{\!\top}\!\! F & -D_u M & D_w & - \tfrac{\gamma(s)^2}{2} I & 0 \\
					\bar{Y}_{1,T}^{d} \left( Z^{d}_{0,T-1} \right)^{\!\top}\!\! F & -B_{\mathrm{cl}}^{d} M & I & 0 & -\tfrac{1}{2} Q
				\end{bmatrix} \prec 0
		\end{equation}}
		\hrulefill 
	\end{table*}
	\!\!\!Apart from \eqref{eq:saturation_LMI} and \eqref{eq:consistency_equation}, the condition that has to be imposed to design a feedback gain $K$ accounting for such an $\ell_{2}$-gain in data-driven control design is dictated by \eqref{eq:l2_gain_LMI}, reported on the next page, leading to the following.
	\begin{thm}\label{thm3}
		 Let $w$ in \eqref{eq:plant_equations} verify \eqref{eq:ell2_bound}. Given the dataset $\mathcal{D}_{T}$ satisfying \eqref{eq:length} and \eqref{eq:rank_condition}, if there exist $Q=Q^{\top}$, a diagonal matrix $M$, $F \in \mathbb{R}^{(n_{u}+n_{x})\times n_{x}}$ and $N \in \mathbb{R}^{n_{u}\times n_{x}}$ satisfying \eqref{eq:saturation_LMI}, \eqref{eq:consistency_equation} and \eqref{eq:l2_gain_LMI}, then the feedback gain $K$ in \eqref{eq:data_driven_K} asymptotically (i.e., for $T \rightarrow \infty$) guarantees $(i)$ closed-loop exponential stability of the origin and that $(ii)$ the $\ell_{2}$-gain bound in \eqref{eq:l2_gain_condition} holds.
	\end{thm}
	\begin{pf}
		Also in this case the Lyapunov candidate in \eqref{eq:lypaunov}, with $P=P^{\top}=Q^{-1}$, is positive definite thanks to the features of the element in $(1,1)$ of \eqref{eq:l2_gain_LMI}. The exponential stability of the origin and property \eqref{eq:l2_gain_condition} can then be established by following the same steps to those carried out in \cite{Seu:23a} (see the proof of Theorem 1), with the only difference that the translation of model-based conditions into their data-driven counterpart is here performed by relying on the asymptotic equivalence of the \textquotedblleft true\textquotedblright \ closed-loop system and the instrument-based one represented by \eqref{eq:data_based_close_loop}. Those steps lead to
 		\begin{equation}\label{eq:lyap_L2}
		 	x \in \mathcal{E}(Q,s) \Rightarrow \Delta V(x,w)+ \frac{1}{\gamma(s)^2} \|z\|_{2}^{2} - \|w\|_{2}^{2} \leq 0,
		 \end{equation}
	 	with $\Delta V(x)$ defined as in the proof of Theorem~\ref{thm2}. Since \eqref{eq:l2_gain_LMI} implies \eqref{eq:reachable_set_LMI}, then Theorem~\ref{thm2} with $\lVert w \rVert_2 \leq s$ and $x(0)=0$ imply $x(t) \in \mathcal{E}(Q,s)$, $\forall t\geq0$. Then, taking sums of the right inequality in \eqref{eq:lyap_L2}, we obtain \eqref{eq:l2_gain_condition} as to be proven. \hfill $\blacksquare$ 
	\end{pf}
	To minimize $\gamma(s)$, by treating $\gamma^2(s)$ as a decision variable, we can cast the following SDP:
	\begin{equation}\label{eq:min_gamma}
		\begin{aligned}
			 & \underset{\gamma^2(s), Q, U, F, N}{\mathrm{minimize}} ~~ \gamma^2(s)                                                  \\
			 & \qquad \quad \mbox{s.t.} ~\quad \eqref{eq:saturation_LMI}, \eqref{eq:consistency_equation}, \eqref{eq:l2_gain_LMI},
		\end{aligned}
	\end{equation}
	which is still tractable when large datasets are used, by virtue of the adopted instrumental variable scheme.
	\section{Numerical example}\label{sec:numerical}
	\begin{figure}[!tb]
		\centering
				\begin{tabular}{c}
			\subfigure[Maximization of the basin of attraction]{
				\begin{tabular}{cc}
					\includegraphics[scale=.25]{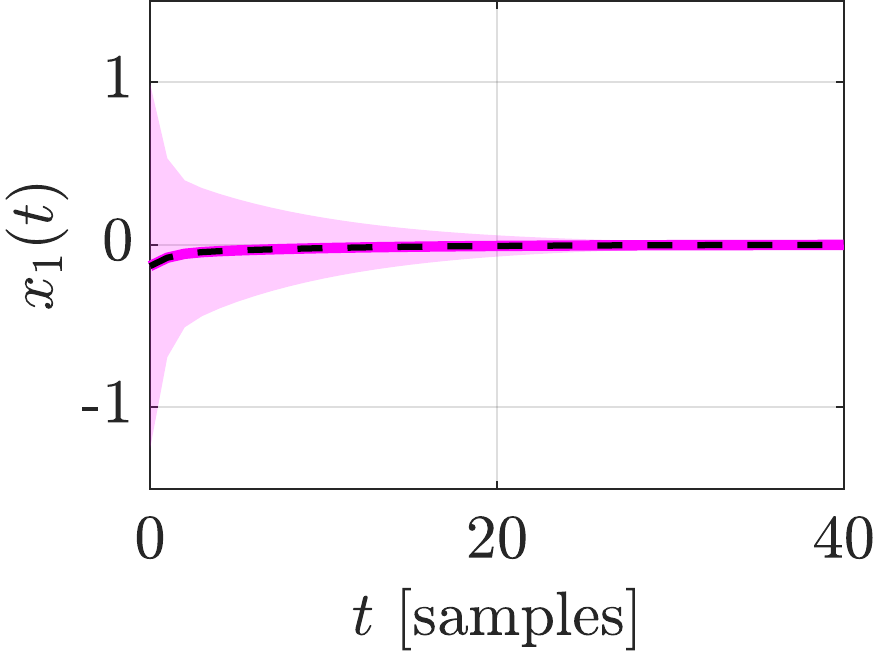} & 
                    \includegraphics[scale=.25]{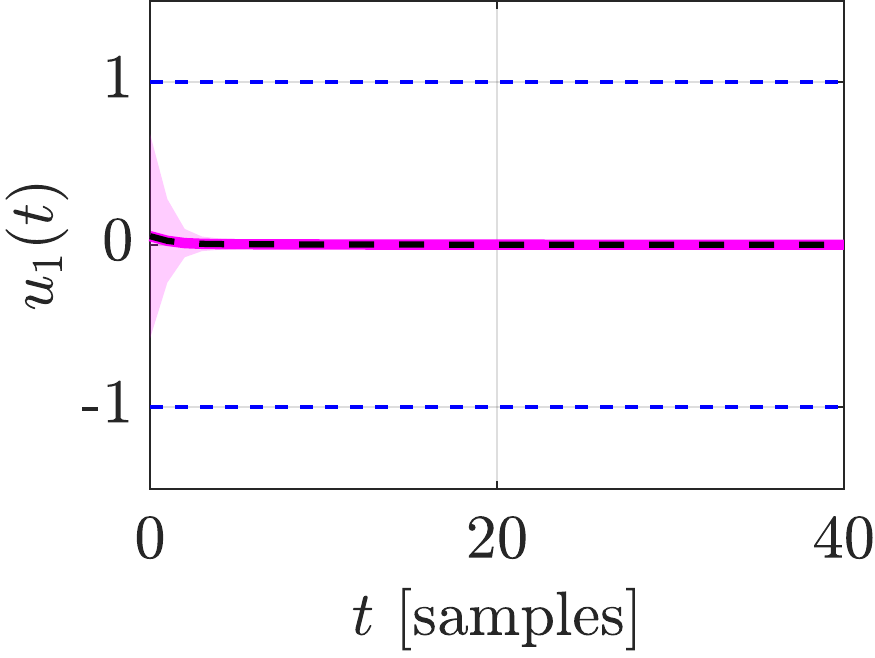}
			\end{tabular}}\\
			\subfigure[Minimization of the estimated reachable set]{
				\begin{tabular}{cc}
					\includegraphics[scale=.25]{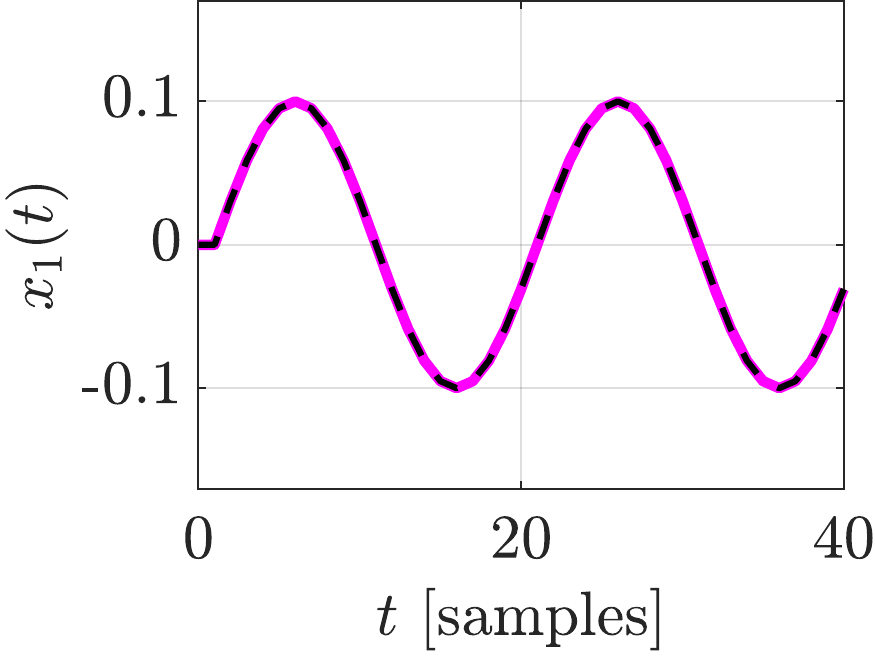} & 
                    \includegraphics[scale=.25]{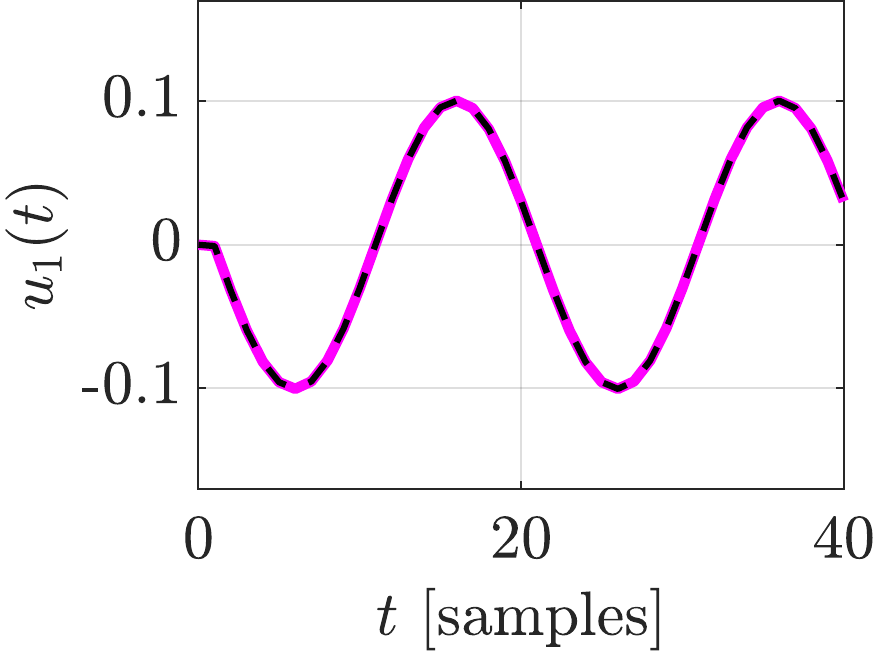}
			\end{tabular}}\\
			\subfigure[Minimization of the $\ell_{2}$-gain $\gamma(s)$]{
				\begin{tabular}{cc}
					\includegraphics[scale=.25]{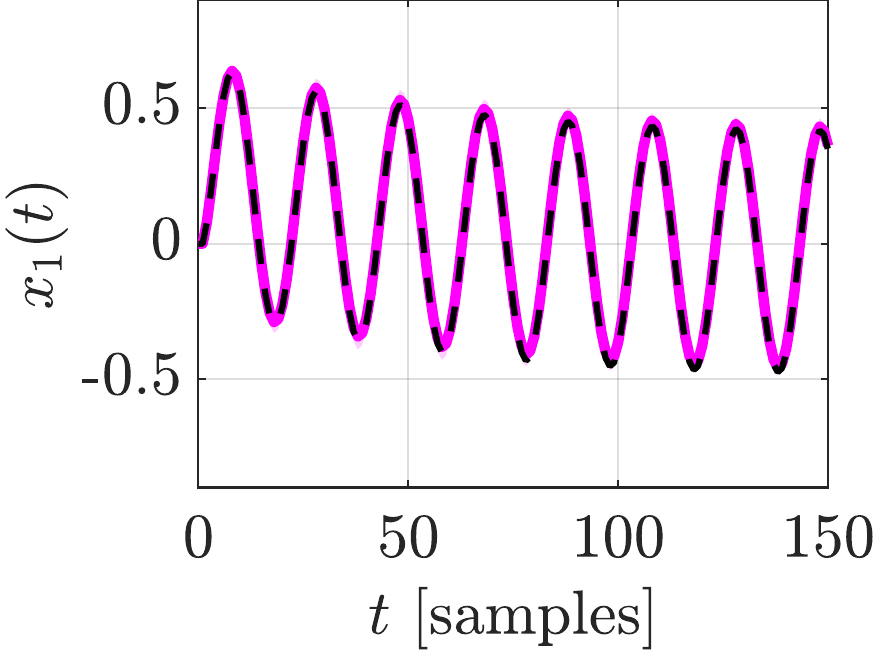} & 
                    \includegraphics[scale=.25]{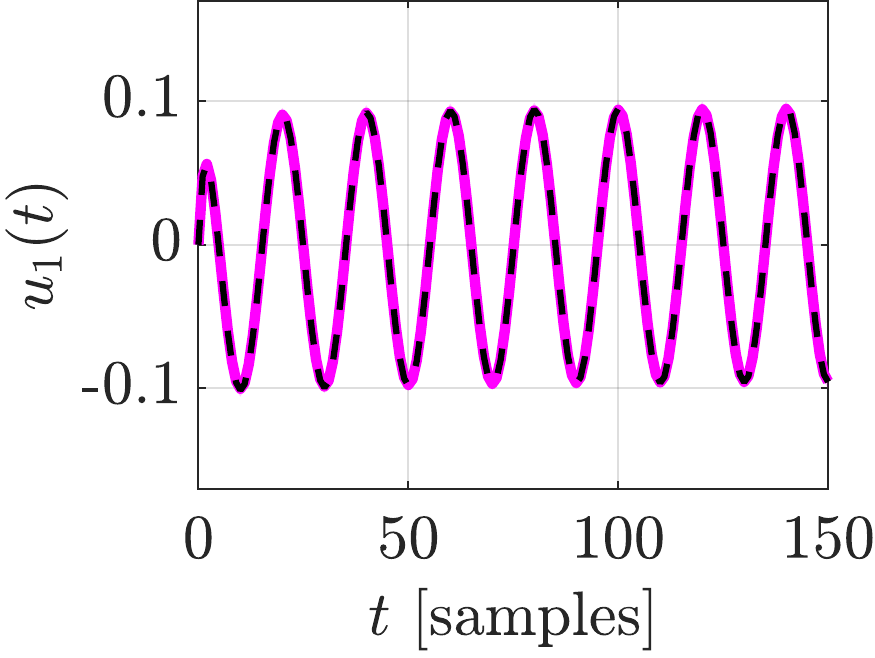}
			\end{tabular}}
		\end{tabular}\vspace{-.2cm}
		\caption{Numerical example: mean and standard deviation (magenta line and shaded area) of the first closed-loop state and input with data-driven controllers over 100 Monte Carlo datasets \emph{vs} behavior with the oracle (black dashed line).}\label{fig:multiple_approaches}
	\end{figure}
	\begin{figure}[!tb]
		\centering
		\includegraphics[scale=.325]{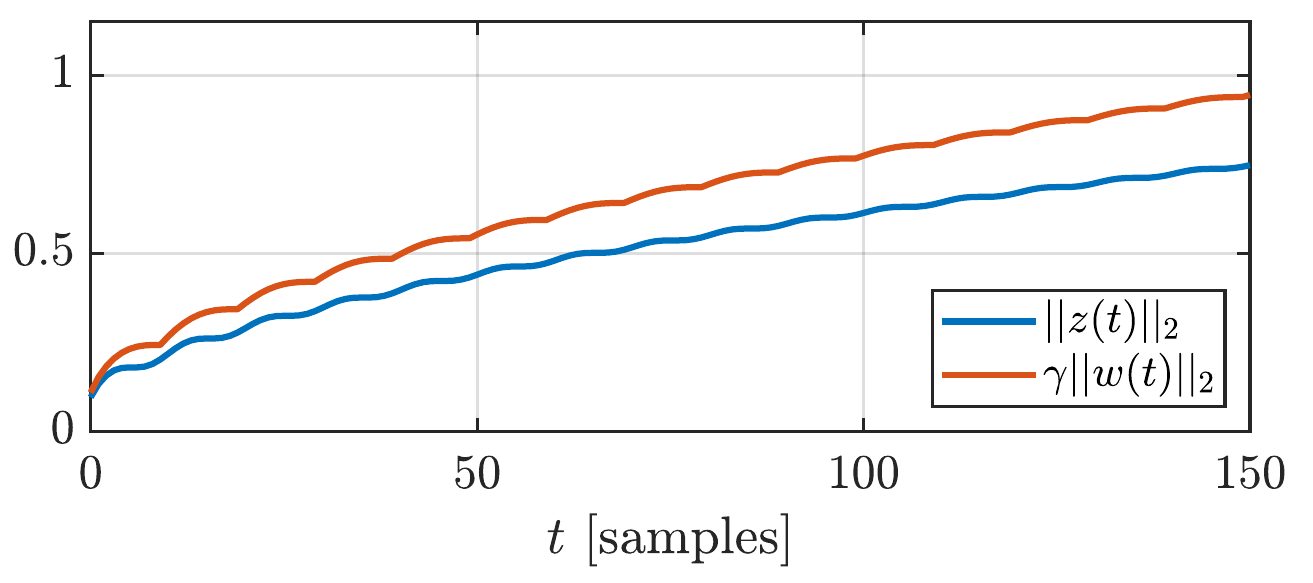}\vspace{-.2cm}
		\caption{Numerical example: mean $\|z\|_{2}$ over $100$ Monte Carlo datasets \emph{vs} $\|w\|_{2}$ scaled.}\label{fig:gain}
	\end{figure}
	Let us consider the open-loop unstable linear system of \cite{Bre:23}, characterized by \eqref{eq:plant_equations} with
		\begin{equation}\label{eq:oracle}
		A = \begin{bmatrix}
			1.01 & 0.01 & 0 \\
			0.01 & 1.01 & 0.01 \\
			0 & 0.01 & 1.01
		\end{bmatrix}, \qquad 
		B = I,
	\end{equation}
	and an input saturation with $\overline{u}_j = 1$, for $j \in \{1,2,3\}$. Moreover, when testing \eqref{eq:trace_min} and \eqref{eq:min_gamma}, we impose
	\begin{equation*}
		C=\begin{bmatrix}
			0 & 1 & 0
		\end{bmatrix},~~D_{u}=\begin{bmatrix}
		-1 & 0 & 1
		\end{bmatrix},~~ D_{w}=\begin{bmatrix}
		0 & 0.3 & -0.8
		\end{bmatrix},
	\end{equation*}
	and set $w(t)$ over closed-loop simulations to
	\begin{equation*}
		w(t)\!=\!0.1\begin{bmatrix}
			\sin\left(\!\frac{t}{10}\pi\!\right) & \sin\left(\!\frac{t}{10}\pi\!+\!\frac{2}{3}\pi\!\right) &
			\sin\left(\!\frac{t}{10}\pi\!+\!\frac{4}{3}\pi\!\right)
		\end{bmatrix}^{\!\top}\!\!.
	\end{equation*}
	We assess the effectiveness of the three data-driven design strategies proposed in this work by collecting $100$ Monte Carlo datasets of length $T = 6000$ in closed-loop\footnote{As in \cite{Bre:23}, we use \eqref{eq:feedback} with $K=I$.} when tracking a uniformly distributed set point in $[-1,1]$. The noise corrupting the output during the data collection phase is zero-mean, Gaussian distributed with covariance $\Sigma_{e}=\sigma_{e}^{2} I$ and $\sigma_{e}=0.1$, yielding a signal-to-noise ratio (SNR) of around $14$~dB. Instead, the controller is tested over noise-free closed-loop experiments, to analyze its effectiveness in realizing the desired control objectives.\\
	Imposing a convergence bound $\eta = 0.995$, we initially evaluate the performance attained when solving \eqref{eq:DoA_optimization} by picking the $100$ initial states drawn randomly from the interval $[-2,2]$. As shown in \figurename{~\ref{fig:multiple_approaches}}\footnote{The shown state and input responses exemplify the behavior of all the other ones.}, the data-driven controller forces state convergence to zero in about 30 time steps for all the tested initial conditions, replicating almost perfectly the mean behavior of the oracle over the tested initial conditions. Similar conclusions on the comparison between our controllers and the oracle can be drawn by looking at the results attained with the controllers designed by solving \eqref{eq:trace_min} and \eqref{eq:min_gamma}, respectively (see \figurename{~\ref{fig:multiple_approaches}}). From \figurename{~\ref{fig:gain}} it is also clear that solving \eqref{eq:trace_min} allows us to satisfy the inequality \eqref{eq:l2_gain_condition}, returning an average value for the $\ell_{2}$-gain $\gamma(s)$ from $w$ to $z$ of $0.86$ over the tests.
	
	\section{Comparison with model identification}\label{sec:result_comparison}
	By focusing on maximizing the origin's basin of attraction, we now compare the performance of the proposed data-driven strategy (DD) with a purely identification-based approach (IB). In this case, the model of the system is the one retrieved by solving the least-squares problem in \eqref{eq:ol_est}, which is used to design the model-based controller based on the certainty equivalence principle.
	
	\textit{Linear system.}
	\begin{figure}[!tb]
		\begin{center}
			\includegraphics[width=0.24\textwidth]{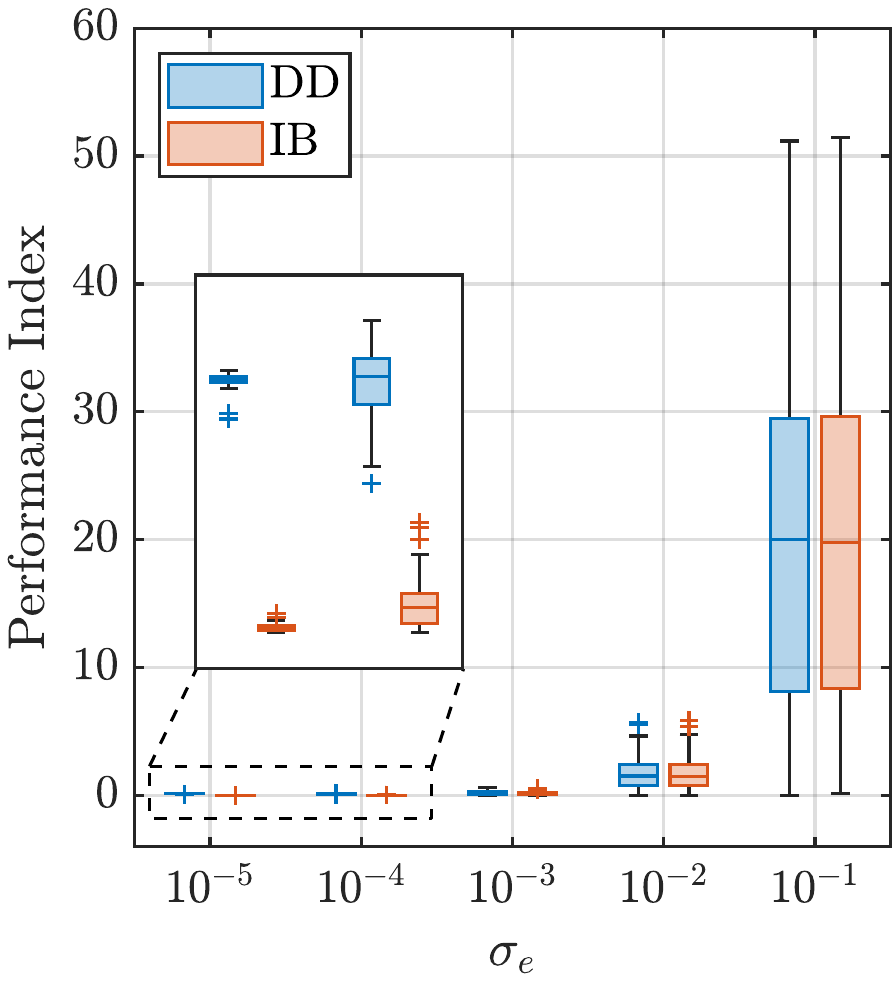}\vspace{-.2cm}
			\includegraphics[width=0.24\textwidth]{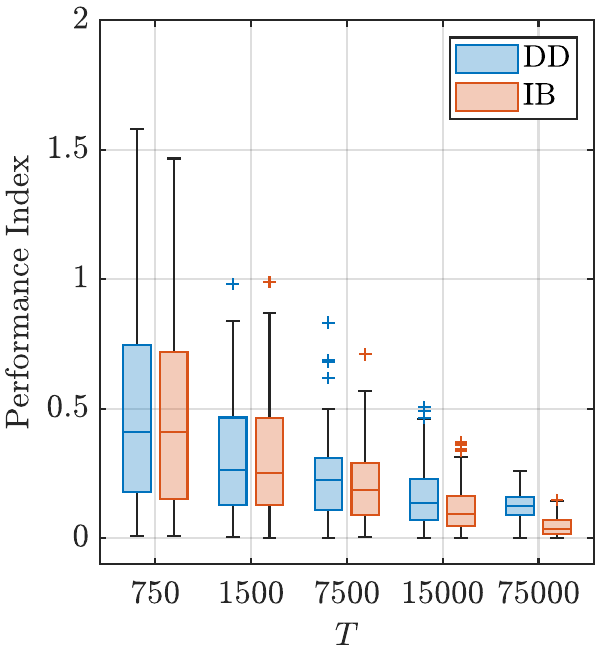} \vspace{-.2cm}
			\caption{Linear system: sensitivity analysis for $\eta = 1$.} 
			\label{fig:eta_1}
		\end{center}
	\end{figure}
		\begin{figure}[!tb]
		\begin{center}
			\includegraphics[width=0.425\textwidth]{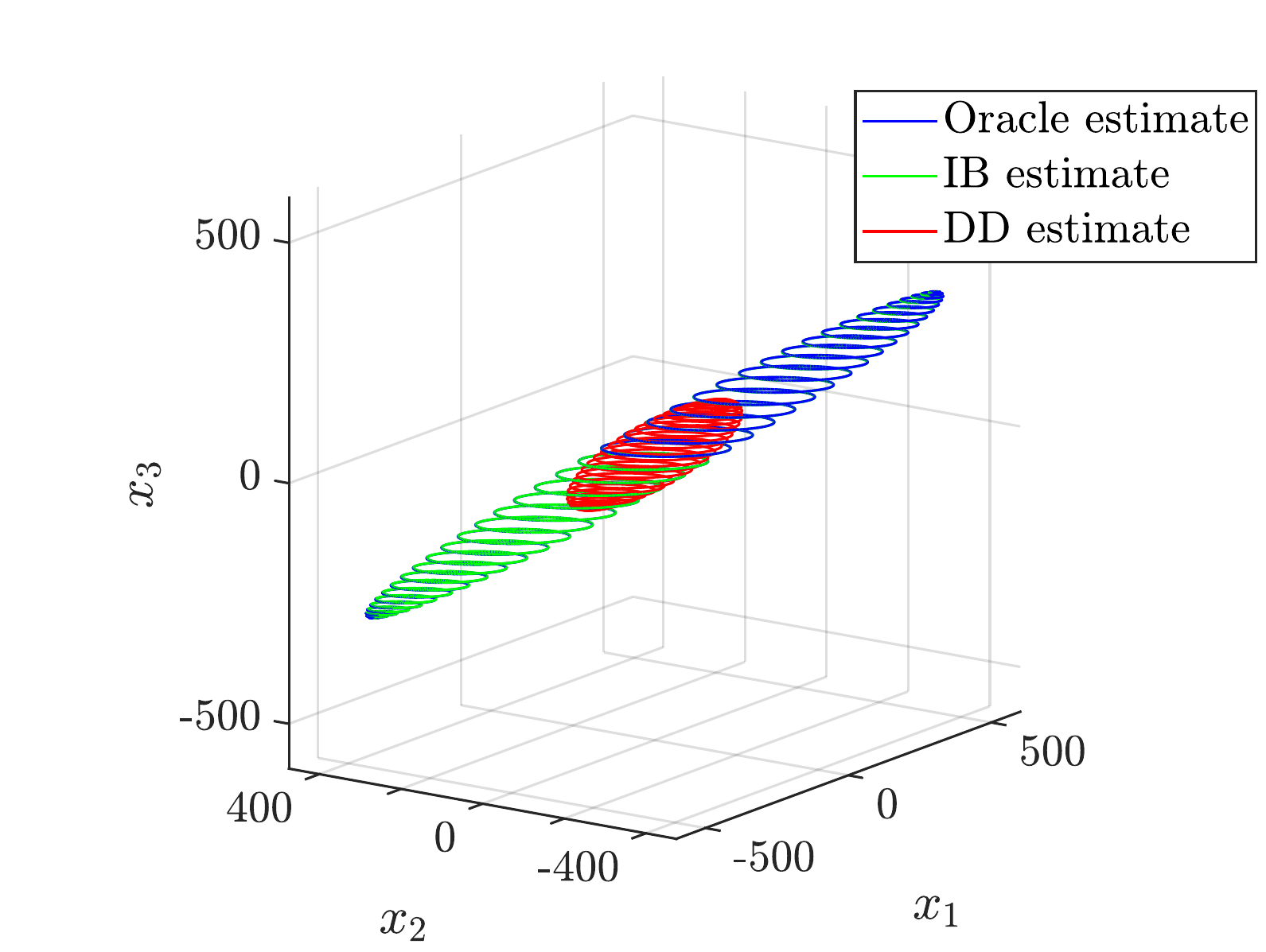}\vspace{-.2cm}
			\caption{Linear system: estimated ellipsoidal basin of attraction for $\sigma_{e}= 10^{-5}$, $T=6000$ and $\eta=1$.} 
			\label{fig:DoA}
		\end{center}
	\end{figure}
		\begin{figure}[!tb]
		\begin{center}
			\includegraphics[width=0.24\textwidth]{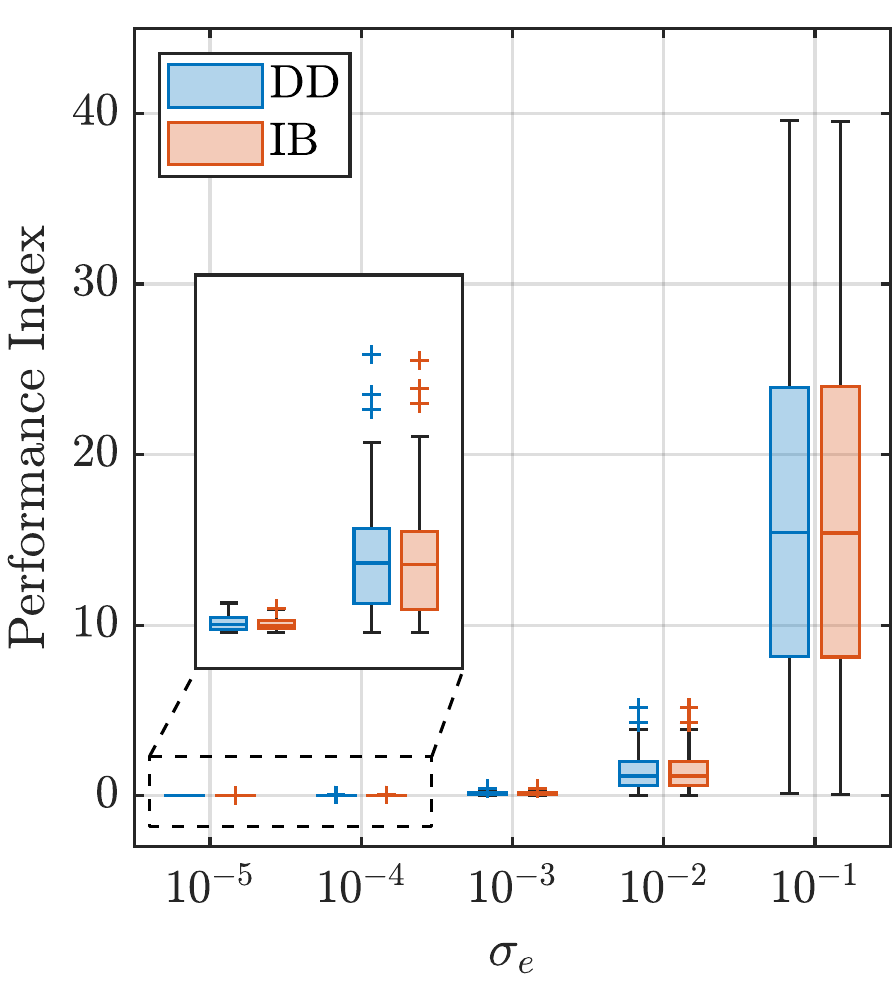}\vspace{-.3cm}
			\includegraphics[width=0.24\textwidth]{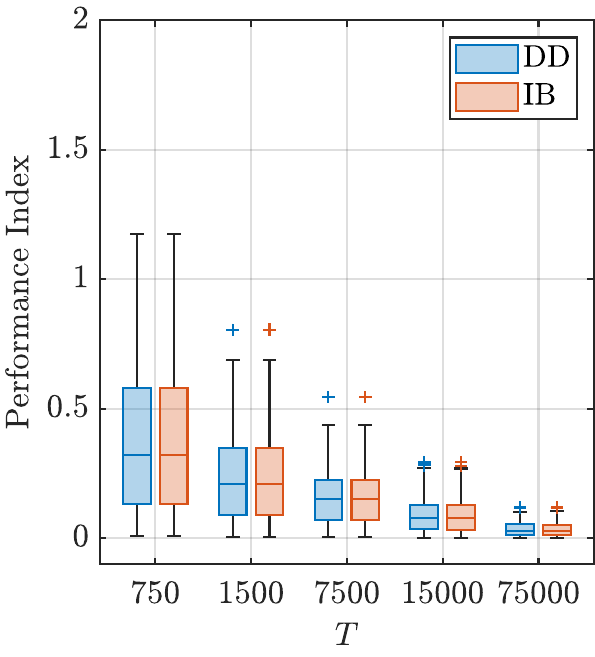} \vspace{-.3cm}
			\caption{Linear system: sensitivity analysis for $\eta = 0.995$.} 
			\label{fig:eta_0995}
		\end{center}
	\end{figure}
	Let us consider the open-loop unstable linear system already introduced in Section~\ref{sec:numerical}, considering the same experimental conditions described before but varying  the length $T$ of the available datasets or changing the noise level. Moreover, let us consider the same performance index used in \cite{Bre:23}, namely
	\begin{equation*}
		\frac{\lvert \alpha^* - \alpha \rvert}{\alpha^*} 100~~~\mathrm{[\%]},
	\end{equation*}
	comparing the solution to \eqref{eq:DoA_optimization} with the \textquotedblleft oracle\textquotedblright \ solution $\alpha^{\star}$, obtained by solving the model-based problem with the true system matrices in \eqref{eq:oracle}.\\
	We firstly impose $\eta=1$ (i.e., not particularly constraining the convergence speed), obtaining the results reported in \figurename{~\ref{fig:eta_1}}. From the latter, it is clear that, while for large noise and small datasets the DD and IB solutions are similar, when the noise is small or the dataset becomes large, the DD strategy returns a solution that does not converge to the oracle one. Such a behavior stems from the SDP solver encountering numerical problems\footnote{The solution of the problem is often inaccurate.} when data directly appears in the LMIs, ultimately leading the DD estimate of the origin's basin of attraction $\mathcal{E}(Q,1)$ to be a strict subset of the one returned by the oracle, which almost coincides to the one returned by the IB strategy (see \figurename{~\ref{fig:DoA}}). Hence, the numerical issues limit the maximum size of the estimated ellipsoid while, as also shown in \figurename{~\ref{fig:DoA}}, they do not affect the ellipsoid orientation. By reducing $\eta$ to $0.995$, thus imposing a faster convergence speed, the numerical issues are nonetheless resolved, leading the DD approach results to be completely equivalent to the IB ones, independently of the noise level and of dataset dimension (see \figurename{~\ref{fig:eta_0995}}).\\
	Comparing the complexity of the SDPs of the two approaches, it is worth to point out that resorting to the DD strategy increases the complexity of the LMI-based problem to be solved. Indeed, $F$ in \eqref{eq:DoA_LMI} is an $(n_x + n_u) \times n_x$ matrix, whereas its dimension reduces to $n_u \times n_x$ when the IB approach is used. All these results indicate that, at least in this case, the model identified by solving \eqref{eq:ol_est} is accurate enough to allow for the design of a controller for the unknown input-saturated system, thus making the shift to the DD approach not particularly advantageous.
	
\begin{figure}[!tb]
	\begin{center}
		\includegraphics[width=0.45\textwidth]{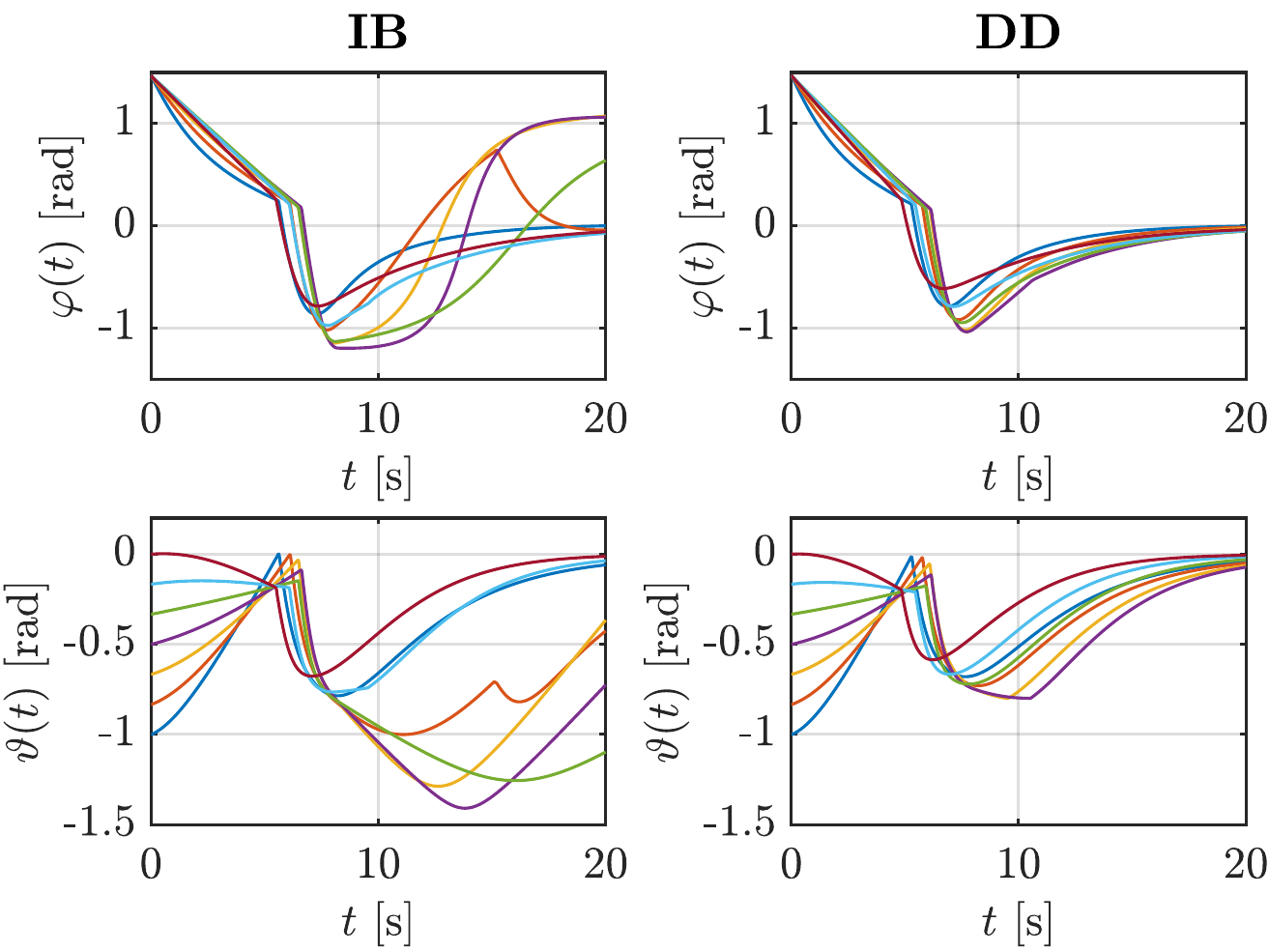}\vspace{-.2cm}
		\caption{Quadcopter control: roll and pitch trajectories for $\varphi(0) = 1.47$ and $\vartheta(0) \in [-1,0]$.} 
		\label{fig:trajectories}
	\end{center}
\end{figure}

	\textit{Nonlinear system.}
		\begin{figure}[!tb]
		\begin{center}
			\includegraphics[width=0.35\textwidth]{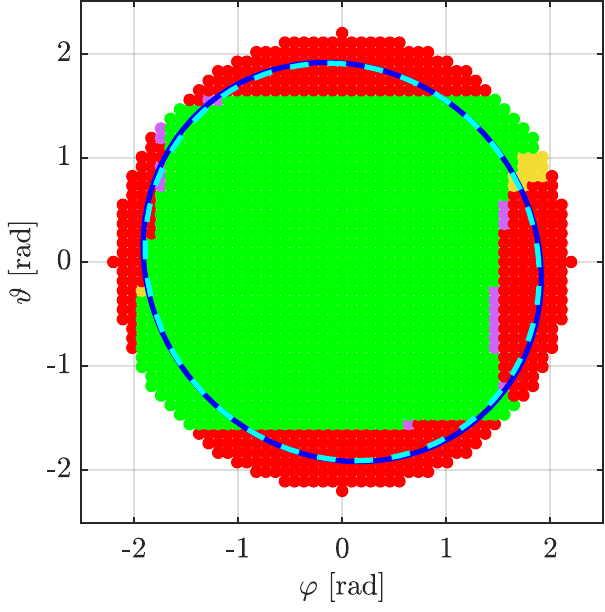}\vspace{-.2cm}
			\caption{Quadcopter control: ellipsoidal estimates of the basin of attraction obtained with the DD (dashed light blue) and IB (solid blue) approaches \emph{vs} initial conditions leading to the satisfaction of \eqref{eq:convergence} at the end of the simulation horizon for both the DD and IB strategies (green dots), for neither of them (red dots), for only the DD approach (purple dots) or for only the IB method (yellow dots).} \label{fig:convergence_experimental}
		\end{center}
	\end{figure}
	For LTI systems, model identification approaches often lead to good closed-loop performance thanks to an accurate reconstruction of the state-space matrices. However, when the true system dynamics are nonlinear, IB approaches devised for linear systems are doomed to make some modelling errors, which often lead to a deterioration in control performance. In this nonlinear setting, we now investigate experimentally whether the proposed DD strategy is able to achieve a better performance by avoiding the full state-space matrix identification. Specifically, we focus on the nonlinear quadrotor system in \cite{For:11}. The approach retains the PID controllers from the mentioned source for controlling the quadrotor position, utilizing both model-based and data-driven versions of \eqref{eq:DoA_optimization} to control the quadrotor's height $z$ [m] and attitude (roll $\varphi$ [rad], pitch $\vartheta$ [rad], and yaw $\psi$ [rad] angles). The objective is to find a state feedback controller \eqref{eq:feedback} that ensures the stability of the quadrotor around a given position while maintaining parallel alignment to the ground.\\
	Within this setting, the data collection phase is carried out by performing closed-loop experiments of $20$ [s] (for a total of $T=2000$ samples) using the controller proposed in \cite{For:11}, by considering step references for the height and yaw angle, random Gaussian references with standard deviation $0.6$ for the roll and pitch angles and corrupting the measured states with a white noise with covariance $\Sigma_{e}=\sigma_{e}I$, with $\sigma_{e}^{2}=10^{-3}$. Meanwhile, \eqref{eq:DoA_optimization} has been augmented with the constraint
	\begin{equation*}
		Q \preceq \kappa^{2} I,
	\end{equation*}  
    with $\kappa^{2}=\pi^{2}/2$, thus accounting for the fact that angles (in absolute value) equal or above $\pi/\sqrt{2}$ imply that the quadrotor is (undesirably) vertical or flipped with respect to the ground. By setting $\eta=0.999$ and considering several initial conditions, we check empirically if the bound in \eqref{eq:convergence} holds for all the attitude angles by considering both the DD and IB approaches at the end of the simulation horizon, while simply checking whether the altitude $z$ at this last instant is positive and lays in an interval between $\pm 1$~[m] from its equilibrium value (to assert that the quadcopter is in a ``non-falling'' mode). As shown in \figurename{~\ref{fig:convergence_experimental}}, differently from the linear case, in this setting the DD approach leads to a wider set of initial conditions verifying \eqref{eq:convergence}. This advantage is further confirmed by \figurename{~\ref{fig:trajectories}}, additionally highlighting that the difference between the DD/IB controllers is due to the inability of the IB solution to steer the quadcopter's attitude to the desired values.

	\section{Conclusions}
		In this paper we have introduced three new data-driven control design strategies for input saturated systems, guaranteeing (asymptotically in the number of data) regional closed-loop stability and either $(i)$ the maximization of the origin's basin of attraction with a  guaranteed convergence rate for the state, or $(ii)$ the minimization of the estimated reachable set for a given bound on uncontrolled (yet measurable) exogenous signals, or $(iii)$ the minimization of the $\ell_{2}$-gain from this uncontrollable signal to the one embedding the desired performance. The presented results show the effectiveness of the proposed strategy when noisy data are employed, showing the potential advantages of a shift to the data-driven rationale in a nonlinear setup. \\
		Future work will be devoted to more complex control schemes, possibly involving an anti-windup compensator.

	\bibliography{ifacconf.bib}             

\begin{thebibliography}{13}
\providecommand{\natexlab}[1]{#1}
\providecommand{\url}[1]{\texttt{#1}}
\providecommand{\urlprefix}{URL }
\expandafter\ifx\csname urlstyle\endcsname\relax
  \providecommand{\doi}[1]{doi:\discretionary{}{}{}#1}\else
  \providecommand{\doi}{doi:\discretionary{}{}{}\begingroup
  \urlstyle{rm}\Url}\fi

\bibitem[{Breschi et~al.(2020)Breschi, Masti, Formentin, and Bemporad}]{Bre:20}
Breschi, V., Masti, D., Formentin, S., and Bemporad, A. (2020).
\newblock {NAW-NET}: neural anti-windup control for saturated nonlinear
  systems.
\newblock In \emph{2020 59th IEEE Conference on Decision and Control},
  3335--3340.

\bibitem[{Breschi et~al.(2023)Breschi, Zaccarian, and Formentin}]{Bre:23}
Breschi, V., Zaccarian, L., and Formentin, S. (2023).
\newblock Data-driven stabilization of input-saturated systems.
\newblock \emph{IEEE Control Systems Letters}, 7, 1640--1645.

\bibitem[{Breschi and Formentin(2020)}]{Bre:20b}
Breschi, V. and Formentin, S. (2020).
\newblock Direct data-driven control with embedded anti-windup compensation.
\newblock In A.M. Bayen, A.~Jadbabaie, G.~Pappas, P.A. Parrilo, B.~Recht,
  C.~Tomlin, and M.~Zeilinger (eds.), \emph{Proceedings of the 2nd Conference
  on Learning for Dynamics and Control}, volume 120 of \emph{Proceedings of
  Machine Learning Research}, 46--54. PMLR.

\bibitem[{da~Silva and Tarbouriech(2005)}]{Sil:05}
da~Silva, J.M.G. and Tarbouriech, S. (2005).
\newblock Antiwindup design with guaranteed regions of stability: an
  {LMI}-based approach.
\newblock \emph{IEEE Transactions on Automatic Control}, 50, 106--111.

\bibitem[{{De Persis} and Tesi(2020)}]{Per:20}
{De Persis}, C. and Tesi, P. (2020).
\newblock Formulas for data-driven control: stabilization, optimality, and
  robustness.
\newblock \emph{IEEE Transactions on Automatic Control}, 65(3), 909--924.

\bibitem[{Dörfler et~al.(2023)Dörfler, Tesi, and {De Persis}}]{Dor:23}
Dörfler, F., Tesi, P., and {De Persis}, C. (2023).
\newblock On the certainty-equivalence approach to direct data-driven {LQR}
  design.
\newblock \emph{IEEE Transactions on Automatic Control}.

\bibitem[{Formentin and Lovera(2011)}]{For:11}
Formentin, S. and Lovera, M. (2011).
\newblock Flatness-based control of a quadrotor helicopter via feedforward
  linearization.
\newblock \emph{IEEE Conference on Decision and Control and European Control
  Conference}, 6171--6176.

\bibitem[{Lasserre(1993)}]{Las:93}
Lasserre, J.B. (1993).
\newblock Reachable, controllable sets and stabilizing control of constrained
  linear systems.
\newblock \emph{Automatica}, 29(2), 531--536.

\bibitem[{Seuret and Tarbouriech(2023{\natexlab{a}})}]{Seu:23a}
Seuret, A. and Tarbouriech, S. (2023{\natexlab{a}}).
\newblock A data-driven approach to the {$L_2$} stabilization of linear systems
  subject to input saturations.
\newblock \emph{IEEE Control Systems Letters}, 7, 1646--1651.

\bibitem[{Seuret and Tarbouriech(2023{\natexlab{b}})}]{Seu:23b}
Seuret, A. and Tarbouriech, S. (2023{\natexlab{b}}).
\newblock Robust data-driven control design for linear systems subject to input
  saturations.
\newblock \emph{arXiv preprint arXiv:2303.04455}.

\bibitem[{Söderström and Stoica(2002)}]{Sod:02}
Söderström, T. and Stoica, P. (2002).
\newblock Instrumental variable methods for system identification.
\newblock \emph{Circuits Systems and Signal Processing}, 21, 1--9.

\bibitem[{Tarbouriech et~al.(2011)Tarbouriech, Garcia, Gomes~da Silva~Jr., and
  Queinnec}]{Tar:11}
Tarbouriech, S., Garcia, G., Gomes~da Silva~Jr., J., and Queinnec, I. (2011).
\newblock \emph{Stability and stabilization of linear systems with saturating
  actuators}.
\newblock Springer-Verlag London Ltd.

\bibitem[{Willems et~al.(2005)Willems, Rapisarda, Markovsky, and {De
  Moor}}]{Wil:05}
Willems, J.C., Rapisarda, P., Markovsky, I., and {De Moor}, B. (2005).
\newblock A note on persistency of excitation.
\newblock \emph{Systems \& Control Letters}, 54(4), 325--329.

\end{thebibliography}
	
\end{document}